\documentclass[a4paper,10pt]{amsart}

\pdfoutput=1
\usepackage{amsmath, amsfonts, amssymb, amsthm, amscd, constants}
\usepackage[pdftex]{graphicx}
\usepackage{color}

\newtheorem{theorem}{Theorem}[section]

\newtheorem{lemma}{Lemma}[section]
\newtheorem{proposition}{Proposition}[section]

\theoremstyle{definition}

\newcommand{\setof}[2]{\left\{#1 \,:\, #2 \right\}}
\newcommand{\given}{\,|\,}
\newcommand{\bgiven}{\bigm|}

\newcommand{\calA}{\mathcal{A}}
\newcommand{\calB}{\mathcal{B}}

\newcommand{\calE}{\mathcal{E}}

\newcommand{\calH}{\mathcal{H}}

\newcommand{\calQ}{\mathcal{Q}}
\newcommand{\calR}{\mathcal{R}}

\newcommand{\calW}{\mathcal{W}}

\newcommand{\bbN}{\mathbb{N}}

\newcommand{\bbZ}{\mathbb{Z}}

\newcommand{\sfZ}{{\sf Z}}

\newcommand{\xmin}{x_{\rm\scriptscriptstyle min}}
\newcommand{\xmax}{x_{\rm\scriptscriptstyle max}}
\newcommand{\ymin}{y_{\rm\scriptscriptstyle min}}
\newcommand{\ymax}{y_{\rm\scriptscriptstyle max}}

\newcommand{\ie}{\textrm{i.e.}}

\newconstantfamily{c}{
symbol=c,
}

\title{Scaling limit of the prudent walk}
\author{V. Beffara, S. Friedli, Y. Velenik}
\email{vbeffara@ens-lyon.fr}
\email{sacha@mat.ufmg.br}
\email{Yvan.Velenik@unige.ch}
\address{UMPA -- ENS Lyon\\46 All\'ee d'Italie\\69364 Lyon Cedex 07, France}
\address{Departamento de Matem\'atica\\UFMG-ICEX C.P. 702\\Belo Horizonte\\30123-970 MG, Brasil}
\address{Section de Math\'ematiques\\Universit\'e de Gen\`eve\\2-4 rue du Li\`evre\\1211 Gen\`eve 4, Suisse}
\begin{document}

\begin{abstract}
We describe the scaling limit of the nearest neighbour prudent walk on $\bbZ^2$, which performs steps uniformly in directions in which it does not see sites already visited. We show that the scaling limit is given by the process $Z_u = \int_0^{3u/7} \bigl( \sigma_1\mathbf{1}_{\{W(s)\geq 0\}} \vec{e}_1 + \sigma_2\mathbf{1}_{\{W(s)< 0\}} \vec{e}_2 \bigr)\mathrm{d} s$, $u\in[0,1]$, where $W$ is the one-dimensional Brownian motion and $\sigma_1,\sigma_2$ two random signs. In particular, the asymptotic speed of the walk is well-defined in the $L^1$-norm and equals $\tfrac37$.
\end{abstract}

\maketitle



\section{Introduction and description of results}
The prudent walk was introduced more than 20 years ago, under the name
\textit{self-directed walk} in~\cite{TD1987b,TD1987a} and \textit{outwardly directed self-avoiding walk} in~\cite{SSK2001},
as a particularly tractable variant of the self-avoiding walk.
Interest in this model has known a
vigorous renewal in recent years, mostly in the combinatorics
community~\cite{D2005,G2006,DG2008,DGGJ2008,S2008, B2008}. The latter works are
concerned with a variant of the original model, more natural from the
combinatorial point of view, obtained by considering the uniform
probability measure on all allowed paths of given length; we'll refer to it as the uniform prudent walk. In the present work, we consider the original (kinetic) model and describe its scaling limit in details.

\medskip
Let $\vec e_1=(1,0)$, $\vec e_2=(0,1)$ denote the canonical basis of $\bbZ^2$. Let us describe the construction of the process associated to the prudent random walk $\gamma_\cdot$. We first set $\gamma_0:=(0,0)$. Assume that we have already constructed\footnote{For any discrete- or continuous-time process $Y_t$, we denote by $Y_{[0,t]} := \{Y_s\}_{0\leq s\leq t}$.} $\gamma_{[0,t]}$, then the distribution of $\gamma_{t+1}$ is obtained as follows. We say that the direction $\vec e\in\{\pm \vec e_1, \pm \vec e_2\}$ is allowed if and only if $\{\gamma_t+k\vec e ,k>0\}\cap \gamma_{[0,t]}=\varnothing$; in other words, allowed directions are those towards which there are no sites that have already been visited. $\gamma_{t+1}$ is chosen uniformly among the neighbours of $\gamma_t$ lying in an allowed direction. Observe that there are always at least two allowed directions.

\medskip
There has been considerable interest in such non-Markovian processes recently. Among the motivations are the necessity of developing specific methods to analyse these processes, and the fact that they sometimes present rather unusual properties. From this point of view, the prudent walk is quite interesting. On the one hand, it is sufficiently simple that a lot of information can be extracted, on the other its scaling limit possesses some remarkable features.

\medskip
\begin{figure}[b]
 \centering
 \pdfximage width \textwidth {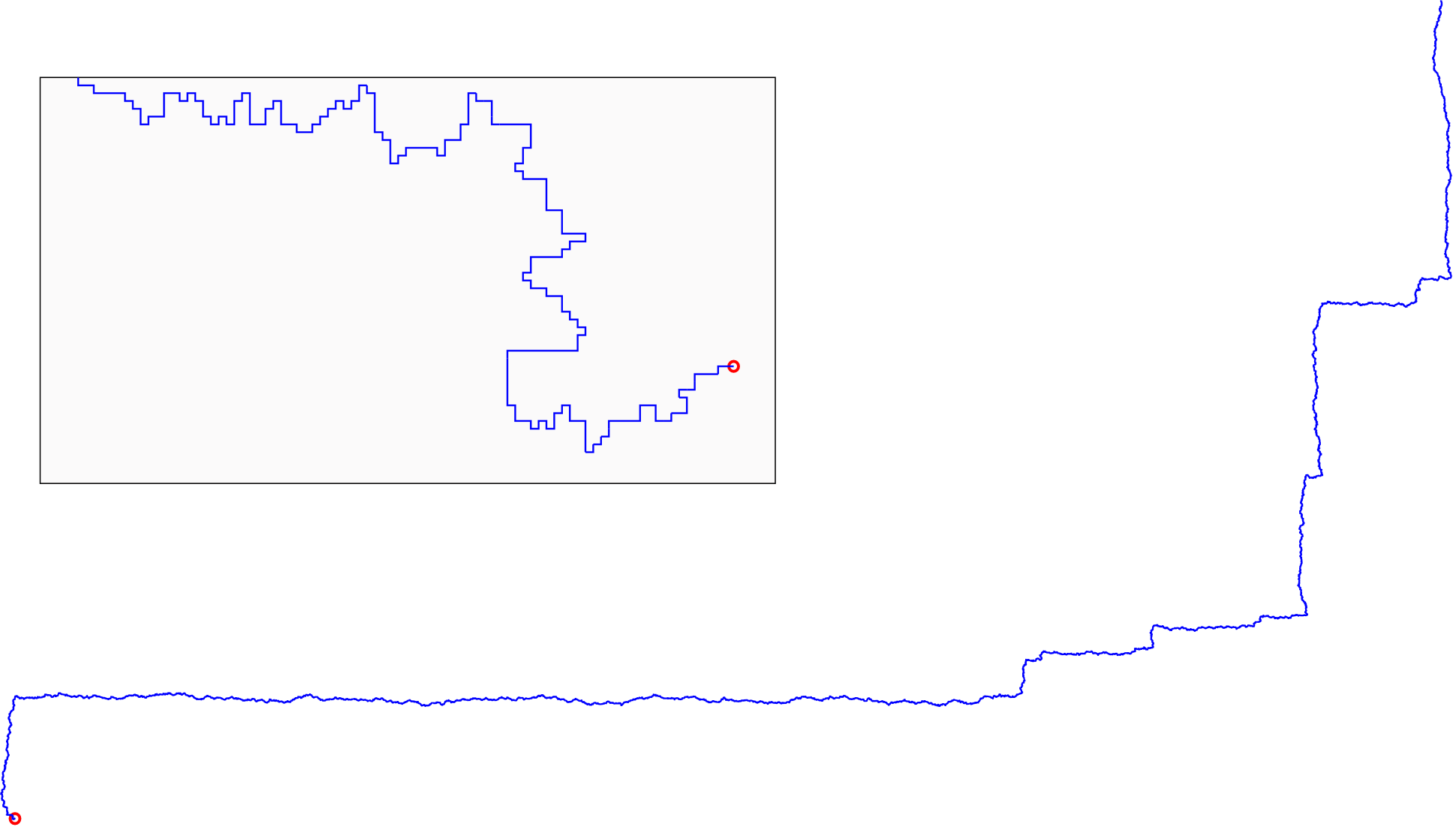}
 \begin{center}
 \pdfrefximage\pdflastximage
\end{center}
\caption{The first $50000$ steps of a prudent walk. The inset is a blow-up of the first few hundred steps near the origin. 
}
 \label{fig:demo}
\end{figure}
Our main result is the following, which describes the scaling limit (in a topology stronger than weak convergence) of the prudent walk.
\pagebreak
\begin{theorem}\label{thm_main}
On a suitably enlarged probability space, we can construct the prudent walk $\gamma_\cdot$, a Brownian motion $W_\cdot$ and a pair of $\pm 1$-valued random variables $\sigma_1,\sigma_2$ such that
\[
\lim_{t\to\infty} P\bigl( \sup_{0\leq s \leq t} \bigl\| \tfrac1t\gamma_s - Z^{\sigma_1,\sigma_2}_{s/t} \bigr\|_{\scriptscriptstyle 2} > \epsilon \bigr) = 0,
\]
where, for $u\in[0,1]$,
\begin{equation}\label{processusZ}
Z^{\sigma_1,\sigma_2}_u := \int_0^{3u/7} \bigl( \sigma_1\mathbf{1}_{\{W(s)\geq 0\}} \vec{e}_1 + \sigma_2\mathbf{1}_{\{W(s)< 0\}} \vec{e}_2 \bigr)\mathrm{d} s.
\end{equation}
The Brownian motion $W_\cdot$ and the two random signs $\sigma_1,\sigma_2$ are independent of each other, and $P(\sigma_1=s,\sigma_2=s') = 1/4$, for $s,s'\in\{-1,1\}$.
\end{theorem}
A typical trajectory of $\gamma_\cdot$ is depicted in~Fig.~\ref{fig:demo}. Let us briefly discuss some of the consequences of Theorem~\ref{thm_main}. To simplify the exposition, we assume that $\sigma_1=\sigma_2=1$ (with no loss of generality, given the obvious symmetries of the prudent walk).
\begin{itemize}
\item Since $\|Z^{1,1}_u\|_{\scriptscriptstyle 1} = 3u/7$, the asymptotic speed of the prudent walk converges to $3/7$ in $L^1$-norm in probability. Note that there is no asymptotic speed in the $L^2$-norm. Actually, using arguments very similar to those of Section~\ref{sec_scalinglimit}, one can prove that the speed converges almost surely to $3/7$.
\item The angle $\alpha_u$ between $Z^{1,1}_u$ and $\vec e_1$ is random, which means that the prudent walk undergoes macroscopic fluctuations in direction, even though it is ballistic. The distribution of $\alpha_u$ can easily be determined, using the arcsine law for Brownian motion:
\begin{align*}
P(\alpha_u \leq x ) &= P(\tan\alpha_u \leq \tan x)\\
&= P\bigl( \frac{\theta^+(W_{[0,u]})}u \geq \frac 1{1+\tan x}\bigr)
= \frac2\pi \arctan \sqrt{\tan x},
\end{align*}
where $\theta^+(W_{[0,u]})$ is the time spent by $W_{[0,u]}$ above $0$.
\item The presence of heavy-tailed random variables (the length of the Brownian excursions) in the limiting process $Z^{1,1}_\cdot$ allows the construction of various natural observables displaying ageing. For example, its first component.
\end{itemize}

\medskip
We emphasise that the scaling limit of the kinetic prudent walk seems to be different from the scaling limit of the uniform prudent walk studied in Combinatorics. Indeed, it is shown in~\cite{B2008}, for a simpler variant of the latter (similar to our corner model, see below) that the scaling limit is a straight line along the diagonal, with a speed in $L^1$-norm approximately given by $0.63$. We expect the same to be true for the real uniform prudent walk.

\medskip
We also observe that our scaling limit is radically different from what is observed for similar continuous random walks that avoid their convex hull, which again are ballistic but, at the macroscopic scale, move along a (random) straight line (\cite{MR2165256}, \cite{MR1961285}).

\bigskip\noindent
\emph{Acknowledgements:} We are grateful to B\'alint T\'oth for suggesting the representation~\eqref{processusZ} for the scaling limit, substantially simpler than the one we originally used. We also thank two anonymous referees for their careful reading. SF was partially supported by CNPq (bolsa de produtividade). Support by a Fonds National Suisse grant is also gratefully acknowledged.

\section{Excursions}
We split the trajectory of the prudent walk into 
a sequence of excursions. 
To a path $\gamma_{[0,t]}$, we associate the bounding rectangle $\calR_t\subset \bbZ^2$ with lower left corner $(\xmin,\ymin)$, 
and upper right corner $(\xmax,\ymax)$. Here,
$\xmin=\min \calA_t$, $\ymin=\min \calB_t$, $\xmax=\max \calA_t$, $\ymax=\max \calB_t$, where
\begin{align*}
\calA_t:=\setof{x\in\bbZ}{\exists y\in\bbZ, (x,y)\in\gamma_{[0,t]}}\,,\,\calB_t:=\setof{y\in\bbZ}{\exists x\in\bbZ, (x,y)\in\gamma_{[0,t]}}\,.
\end{align*}
Notice that $\gamma_t$ always lies on the boundary of its bounding rectangle.
The height and width of $\calR_t$ (in units of lattice sites) are defined by $\calH_t:=|\calA_t|$ and $\calW_t:=|\calB_t|$; in particular $\calH_0=\calW_0=1$.
We say that the path visits a corner at time $s\in\{0,\dots,t\}$ 
if $\gamma_s$ coincides with one of the four corners of $\calR_s$.

\medskip
Without loss of generality, we assume that the first step of the path is in the direction $\vec e_1$, \ie, $\gamma_1=\vec e_1$, and 
define the following sequence of random times: 
\begin{align*}
T_0&:=0\,,\\
U_0&:=\inf\setof{t>0}{\calH_t>{1}}-1\,,
\end{align*}
and for $k\geq 0$,
\begin{align*}
T_{k+1} &= \inf\setof{t>U_{k}}{\calW_t>\calW_{t-1}}-1\,,\\
U_{k+1} &= \inf\setof{t>T_{k+1}}{\calH_t>\calH_{t-1}}-1\,.
\end{align*}
During the time interval $(T_k,U_k]$, the walk makes an excursion, 
denoted $\calE_k^{\rm\scriptscriptstyle v}$, along one of the vertical sides of $\calR_t$.
During the time interval $(U_k,T_{k+1}]$, the walk makes an excursion, denoted
$\calE_k^{\rm\scriptscriptstyle h}$,
along one of the horizontal sides of $\calR_t$.
Two relevant quantities are
the horizontal displacement of $\calE_k^{\rm\scriptscriptstyle v}$: $X_{k} := \calW_{T_{k+1}} - \calW_{T_k}$, and
the vertical displacement of $\calE_k^{\rm\scriptscriptstyle h}$: $Y_k := \calH_{T_{k+1}} - \calH_{T_k}$.
Observe that by construction, $X_k\geq 1$ and $Y_k\geq 1$.

\begin{figure}
\includegraphics[width=8cm]{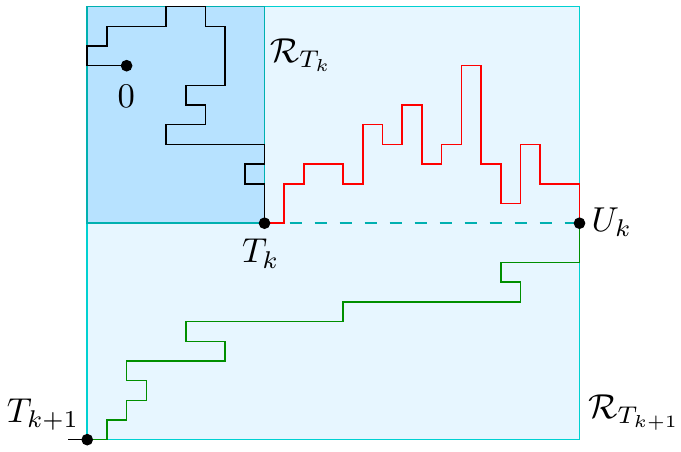}
\caption{The excursions and their associated times.}
\label{Figure1}
\end{figure}

\subsection{The effective random walk for the excursions}
\label{ssec-effective}
As Figure \ref{Figure1} suggests,
the study of excursions along the sides of the rectangle can be reduced to that
of the excursions of an effective one-dimensional random walk, with geometric increments.

\medskip
Let $\xi_1,\xi_2,\dots$ be an i.i.d. sequence, $\xi_i\in\bbZ$, with 
$P(\xi_i=k)=\tfrac13(\tfrac12)^{|k|}$.
Let $S_0=0$, $S_n=\xi_1+\dots+\xi_n$. 
We call $S_n$ the effective random walk.
Let $\eta_L$ be the first exit time of $S_n$ from the interval $[0,L-1]$: $\eta_L:=\inf\{n>0:S_n\not\in [0,L-1]\}$.
The following lemma will be used repeatedly in the paper.
\begin{lemma}\label{lem_effective} For all $k\geq 0$, $m\geq 1$,
\begin{align}
P(X_k=m|\gamma_{[0,T_k]})&=P(\eta_{\calH_{T_k}}=m)\,,\label{E1}\\
P(Y_k=m|\gamma_{[0,U_k]})&=P(\eta_{\calW_{U_k}}=m)\,.
\end{align}
\end{lemma}

\begin{figure}
\scalebox{1.2}{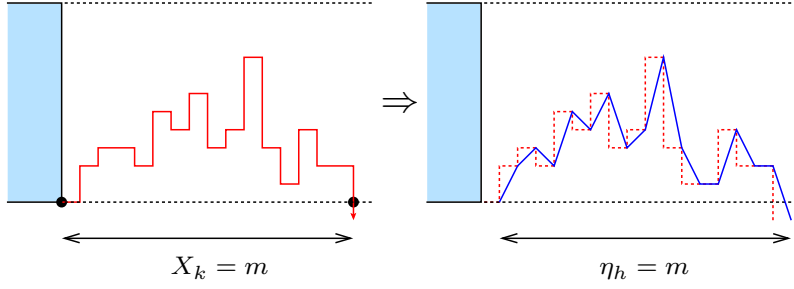}
\caption{The decomposition of the excursion $\calE_k^{\rm\scriptscriptstyle v}$.
The full segments on the right image represent the increments of the effective random walk $S_n$.
Observe that $S_n$ starts at $\gamma_{T_k+1}$, and that it exits $[0,h-1]$ at time $m$ by jumping to any point on the negative axis.}
\label{Fig2}
\end{figure}

\begin{proof} We show \eqref{E1}. 
Assume that the walk is at the lower right corner of the rectangle at time $T_k$.
On $\{X_k=m\}\cap\{\calH_{T_k}=h\}$, the excursion $\calE_k^{\rm\scriptscriptstyle v}$ can be decomposed into
$(\alpha_1,\dots,\alpha_{m-1},\beta)$, where each $\alpha_i$ is an elementary increment 
made of a vertical segment of length $l_i\in[-h+1,h-1]\cap \bbZ$ and of a 
horizontal segment of length $1$ pointing to the right, and $\beta$ is a  
purely vertical segment of length $l\in[-h,+h]\cap \bbZ$, which ensures that the excursion reaches either the top or the bottom of the rectangle by time $m$ (see Figure \ref{Fig2}). 
We thus have
$$
P(X_k=m|\gamma_{[0,T_k]})=\sum_{(\alpha_1,\dots,\alpha_{m-1},\beta)}
p(\alpha_1)\dots p(\alpha_{m-1})p(\beta)\,,
$$
where the sum is over all possible sets of such elementary increments, and $p(\alpha_i)=\tfrac13(\tfrac12)^{|l_i|}$, $p(\beta)=\tfrac13(\tfrac12)^{|l|-1}$. Therefore, in terms of the effective random walk with increments $\xi_i$,
$$
p(\alpha_1)\dots p(\alpha_{m-1})=P(\xi_1=l_1,\dots,\xi_{m-1}=l_{m-1})\,.
$$
Moreover, if $l<0$ (as on Figure \ref{Fig2}), then
$$
p(\beta)=\tfrac13(\tfrac12)^{|l|}
=\tfrac13\sum_{j\geq |l|+1}(\tfrac12)^j\equiv P(S_{m}<0\given S_{m-1}=|l|)\,.$$
Similarly, if $l\geq 0$ 
$$
p(\beta)= P(S_{m}\geq h\given S_{m-1}=h-l-1)\,.
$$
This shows that 
$$P(X_k=m|\gamma_{[0,T_k]})=P(S_{1}\in [0,h-1],\dots,S_{m-1}\in [0,h-1],S_m\not\in [0,h-1])\,,$$
which is \eqref{E1}.
\end{proof}

\section{Crossings and visits to the corners}\label{SectionCrossings}
In this section, we show that, almost surely, the prudent walk eventually
visits exactly one of the corners of its bounding rectangle an infinite number of times. This excludes, in particular, 
the possibility of winding around the origin an infinite number of times, and will
allow in the sequel 
to restrict the study of the prudent walker to a single corner of the rectangle.

\medskip
Let $A_k$ denote the event in which the 
$k$th excursion crosses at least one side of $\calR_{T_k}$, that is if the walk is at a bottom (resp. top) corner of the box at time $T_k$ and at a top (resp. bottom) corner of the box at time $U_{k}$, or if it is at a right (resp. left) corner of the box at time $U_k$ and at a left (resp. right) corner of the box at time $T_{k+1}$.
\begin{proposition}\label{Prop1}
There exists a constant $\Cl[c]{const1}>0$ such that
\begin{equation}\label{E10}
P(A_k)\leq \frac{\Cr{const1}}{k^{4/3}}\,
\end{equation}
for all large enough $k$. In particular, $P(\limsup_k A_k)=0$, and a.s.\ exactly one of the corners is visited infinitely many times.
\end{proposition}
For $k$ large, the bounding rectangle has long sides, which makes crossing unlikely. 
As is well known from the gambler's ruin estimate for the simple symmetric random walk, the probability of first leaving an interval of length $L$ at the opposite end is of order $L^{-1}$. To obtain \eqref{E10}, it will therefore be sufficient 
to show that the sides of the rectangle at the $k$th visit to a corner
grow superlinearly in $k$. 

\medskip
Before this, we give a preliminary 
result for the effective random walk $S_n$ starting at $0$. Let as before
$\eta_L$ denote the first exit time of $S_n$ from the interval $[0,L-1]$. We will also use
$\eta_\infty:=\inf\{n>0:S_n< 0\}$. 
\begin{lemma}\label{L1}
There exists a constant $\Cl[c]{const2}>0$ such that for all $L\geq 1$,
\begin{equation}\label{E9}
P(\eta_L\geq n)\geq \frac{\Cr{const2}}{\sqrt{n}}\quad\text{for all integer }n\leq \Cr{const2}L^{4/3}\,.
\end{equation}
\end{lemma}
\begin{proof}
Define $\eta_L^{\rightarrow}:=\inf\{n>0:S_n\geq L-1\}$, so that 
\begin{align*}
 P(\eta_L\geq n)&=P(\eta_\infty \geq n,\eta_L^{\rightarrow}\geq n)\\
&=P(\eta_\infty \geq n)-P(\eta_L^{\rightarrow}< n,\eta_\infty \geq n)\\
&\geq P(\eta_\infty \geq n)-P(\eta_L^{\rightarrow}< n)\,.
\end{align*}
By the gambler's ruin estimate of Theorem 5.1.7. in \cite{LL2008},
\begin{equation}\label{GREst1}
P(\eta_\infty \geq n)\geq \frac{\Cl[c]{const3}}{\sqrt{n}}
\end{equation}
for some constant $\Cr{const3}>0$. On the other hand, 
$$P(\eta_L^{\rightarrow}< n)
\leq P(\max_{1\leq j\leq n} |S_j|\geq L)\leq \frac{E[|S_n|^2]}{L^2}=\frac{2 n}{L^2}
$$
The second inequality follows from the Doob-Kolmogorov Inequality.
\end{proof}

As a consequence, we can show that by the $k$th visit to the corner, the sides of the rectangle have grown by at least $k^{4/3}$. A refinement of the method below actually shows that $\calW_{T_k}$ and $\calH_{T_k}$ grow like $k^2$.

\begin{lemma}\label{L2}
There exist positive constants $\Cl[c]{const4}$ and $\Cl[c]{const2p}$ such that
\begin{align}
P(\calW_{T_k}< \Cr{const2p}k^{4/3})&\leq \exp\big(-\Cr{const4} k^{1/3}\big)\,,\\
P(\calH_{T_k}< \Cr{const2p}k^{4/3})&\leq \exp\big(-\Cr{const4} k^{1/3}\big)\,.
\end{align}
\end{lemma}

\begin{proof} Let $m=\lfloor k/2\rfloor$.
Since $X_i\geq 1$ and $Y_i\geq 1$ for all $i$, we have that $\calH_{T_j}\geq m$, $\calW_{T_j}\geq m$, for all $j\geq m$.
We consider the width of the rectangle at time $T_k$.
Let $I_j$ denote the indicator of the event ${\{X_j\geq \Cr{const2}m^{4/3}\}}$.
We have, for all $j\geq m$,
\begin{align*}
P(I_j=1\given \gamma_{[0,T_j]})&=P(I_j=1\given \calH_{T_j})\\
&=P(\eta_{\calH_{T_j}}\geq \Cr{const2}m^{4/3})\\
&\geq P(\eta_{m}\geq \Cr{const2}m^{4/3})\\
\text{by }\eqref{E9}&\geq \Cr{const2}^{-3/2}m^{-2/3}\equiv p\,.
\end{align*}
Therefore, the $I_j$s can be coupled 
to Bernoulli variables of parameter $p$. 
Since  $\calW_{T_k}=\sum_{j=0}^{k-1}X_j\geq \Cr{const2}m^{4/3} \sum_{j=m}^{k-1}I_j$, we get
\begin{align*}
P(\calW_{T_k}< 2^{-4/3}\Cr{const2}k^{4/3})&\leq P(I_j=0\,,\forall j=m,\dots,k-1)\\
&\leq(1-p)^{m}\\
&\leq e^{-\Cr{const4}k^{1/3}}\,.
\end{align*}
\end{proof}

\begin{proof}[Proof of Proposition \ref{Prop1}.]
Consider the event $A_k^{\rm\scriptscriptstyle v}$ in which $\calE_k^{\rm\scriptscriptstyle v}$
crosses the (vertical) side of the rectangle.  We have
\begin{align*}
P(A_k^{\rm\scriptscriptstyle v})
&\leq P\bigl(\calH_{T_k}< \Cr{const2p} k^{4/3}\bigr)+P\bigl(A_k^{\rm\scriptscriptstyle v}\bgiven
\calH_{T_k}\geq \Cr{const2p} k^{4/3}\bigr)
\end{align*}
The first term is treated with Lemma \ref{L2}.
Proceeding as in the proof of Lemma~\ref{lem_effective}, we see that, in terms of the effective random walk $S_n$, 
the second term is the probability that $S_n$ (started at zero) reaches $[\Cr{const2p}k^{4/3},+\infty)$ before becoming negative.
Therefore, again by Theorem 5.1.7. in \cite{LL2008}, 
there exists a constant $\Cl[c]{consttt}>0$ such that 
$$P\bigl(A_k^{\rm\scriptscriptstyle v}\bgiven\calH_{T_k}\geq 
\Cr{const2p} k^{4/3}\bigr)\leq \frac{\Cr{consttt}}{k^{4/3}}\,.$$
This shows \eqref{E10}. The second claim follows by the Borel-Cantelli Lemma and by recurrence of the effective walk.
\end{proof}

\begin{figure}
\scalebox{0.8}{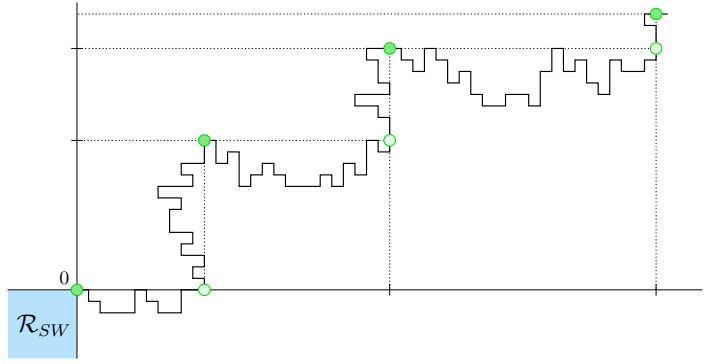}
\caption{A trajectory of the corner process $\widehat\gamma_\cdot$.}
\label{Fig3}
\end{figure}

Proposition \ref{Prop1} implies that the prudent walker behaves qualitatively, asymptotically, in the same way as a simplified model in which the evolution is along the corner of the infinite 
rectangle $\calR_{SW}:=\{(x,y)\in \bbZ^2:x\leq 0,y\leq 0 \}$. 
Namely, redefine a direction 
$\vec e\in\{\pm \vec e_1, \pm \vec e_2\}$ to be allowed if and only if 
$\{\gamma_t+k\vec e ,k>0\}\cap (\gamma_{[0,t]}\cup \calR_{SW})=\varnothing$. The paths of the corner process $\widehat{\gamma}_\cdot$ are obtained by choosing at each step a nearest 
neighbour, uniformly in the allowed directions. 
A trajectory of the corner model is depicted on Figure~\ref{Fig3}.
As before,
$\widehat{\gamma}_\cdot$ can be decomposed into a 
concatenation of excursions along the vertical and horizontal sides of $\calR_{SW}$, denoted respectively 
${\widehat{\calE}}^{\rm\scriptscriptstyle v}_k$, ${\widehat{\calE}}^{\rm\scriptscriptstyle h}_k$. 
Since there is no crossing possible here, these are independent.

\begin{figure}[t]
\begin{center} 
\scalebox{0.33}{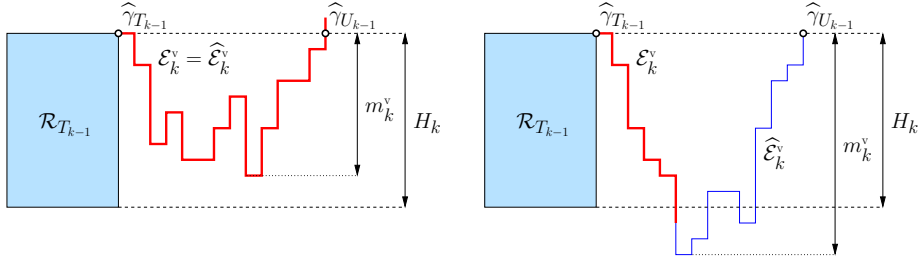}
\end{center}
\caption{The construction of the vertical excursion $\calE_k^{\rm\scriptscriptstyle v}$ of $\gamma_\cdot$ from the vertical excursion $\widehat{\calE}_k^{\rm\scriptscriptstyle v}$ of $\widehat{\gamma}_\cdot$. When $\widehat{\calE}_k^{\rm\scriptscriptstyle v}$ is compatible with $\calR_{T_{k-1}}$, it is kept (left); when it is not compatible, it is truncated (right).}
\label{fig:coupling}
\end{figure}

\medskip
We relate the true excursions of the prudent walk to those of the corner model, by means of the following coupling.
As before, we assume that the first step is horizontal: $\widehat{\gamma}_1=(1,0)$.
Consider a given realisation of 
$({\widehat{\calE}}^{\rm\scriptscriptstyle v}_k$, ${\widehat{\calE}}^{\rm\scriptscriptstyle h}_k)_{k\geq 1}$, from which we construct $({{\calE}}^{\rm\scriptscriptstyle v}_k, {{\calE}}^{\rm\scriptscriptstyle h}_k)_{k\geq 1}$, with the same distribution as the excursions of the prudent walk. Let (see Fig.~\ref{fig:coupling})
\[
m_k^{\rm\scriptscriptstyle v} := \widehat{\gamma}_{T_{k-1}}(2) - \min_{T_{k-1} \leq t \leq U_{k-1}} \widehat{\gamma}_t(2), \quad
m_k^{\rm\scriptscriptstyle h} := \widehat{\gamma}_{U_{k-1}}(1) - \min_{U_{k-1} \leq t \leq T_k} \widehat{\gamma}_t(1),
\]
where we have denoted by $x(1)$ and $x(2)$ the coordinates of $x\in\bbZ^2$.
If $m_k^{\rm\scriptscriptstyle v}> \ell$, let $t_\ell$ denote the first time $T_{k-1} \leq t \leq U_{k-1}$ such that $\widehat{\gamma}_{T_{k-1}}(2)-\widehat{\gamma}_t(2)=\ell+1$, and let 
$\mathrm{Trunc}_\ell(\widehat{\calE}_k^{\rm\scriptscriptstyle v})$ denote the restriction of 
$\widehat{\calE}_k^{\rm\scriptscriptstyle v}$ up to time $t_\ell$. Similarly,
one defines the corresponding objects in the horizontal case.
Let $H_0:=0$ and set 
$$
\calE_{1}^{\rm\scriptscriptstyle v}:=
\begin{cases}
\widehat{\calE}_{1}^{\rm\scriptscriptstyle v}&\text{ if }
m_{1}^{\rm\scriptscriptstyle v} \leq H_0\,,\\
\mathrm{Trunc}_{H_0}(\widehat{\calE}_{1}^{\rm\scriptscriptstyle v})&\text{ if }
m_{1}^{\rm\scriptscriptstyle v} > H_0\,.
\end{cases}
$$
Let $X_0$ be the horizontal displacement associated to $\calE_{1}^{\rm\scriptscriptstyle v}$ as in Figure~\ref{Fig2}, and set $W_0=X_0$. We then set
$$
\calE_{1}^{\rm\scriptscriptstyle h}:=
\begin{cases}
\widehat{\calE}_{1}^{\rm\scriptscriptstyle h}&\text{ if }
m_{1}^{\rm\scriptscriptstyle h} \leq W_0\,,\\
\mathrm{Trunc}_{W_0}(\widehat{\calE}_{1}^{\rm\scriptscriptstyle h})&\text{ if }
m_{1}^{\rm\scriptscriptstyle h} > W_0\,.
\end{cases}
$$
Let $k\geq 1$. Assume $({{\calE}}^{\rm\scriptscriptstyle v}_j, {{\calE}}^{\rm\scriptscriptstyle h}_j)_{j=1,\dots,k}$ have already been defined. To each vertical excursion ${\calE}^{\rm\scriptscriptstyle v}_j$, we associate a horizontal displacement $X_j$ as in Figure~\ref{Fig2}. Similarly, to each horizontal excursion ${\calE}^{\rm\scriptscriptstyle h}_j$, we associate a vertical displacement $Y_j$.

Set $H_k:=Y_0+\dots+Y_{k-1}$, and
$$
\calE_{k+1}^{\rm\scriptscriptstyle v}:=
\begin{cases}
\widehat{\calE}_{k+1}^{\rm\scriptscriptstyle v}&\text{ if }
m_{k}^{\rm\scriptscriptstyle v} \leq H_k\,,\\
\mathrm{Trunc}_{H_k}(\widehat{\calE}_{k+1}^{\rm\scriptscriptstyle v})&\text{ if }
m_{k}^{\rm\scriptscriptstyle v} > H_k\,.
\end{cases}
$$
Let now $W_k:=X_0+\dots+X_k$, where $X_k$ is the horizontal displacement of the vertical excursion $\calE_{k+1}^{\rm\scriptscriptstyle v}$, and set
$$
\calE_{k+1}^{\rm\scriptscriptstyle h}:=
\begin{cases}
\widehat{\calE}_{k+1}^{\rm\scriptscriptstyle h}&\text{ if }
m_{k}^{\rm\scriptscriptstyle h} \leq W_k\,,\\
\mathrm{Trunc}_{W_k}(\widehat{\calE}_{k+1}^{\rm\scriptscriptstyle h})&\text{ if }
m_{k}^{\rm\scriptscriptstyle h} > W_k\,.
\end{cases}
$$
By construction, vertical and horizontal excursions of the prudent walk alternate.
Therefore, $({{\calE}}^{\rm\scriptscriptstyle v}_k, {{\calE}}^{\rm\scriptscriptstyle h}_k)_{k\geq 1}$ is 
distributed exactly as the excursions of the true prudent walk up to the appropriate reflection.
By Proposition \ref{Prop1}, almost surely, $\calE_{k}^{\rm\scriptscriptstyle h}=\widehat{\calE}_{k}^{\rm\scriptscriptstyle h}$
and $\calE_{k}^{\rm\scriptscriptstyle v}=\widehat{\calE}_{k}^{\rm\scriptscriptstyle v}$ for all but a finite number of $k$s.

\medskip
According to Proposition~\ref{Prop1}, the prudent walker $\gamma_\cdot$ eventually fixates in one of the four quadrants and couples with the trajectory of the corner model $\hat \gamma_\cdot$. Obvious lattice symmetries 
present in the model imply that the quadrant is chosen with uniform probability $1/4$. In order to lighten notations, we shall assume that the chosen quadrant is the first one.
\begin{theorem}\label{thm_cornervsprudent}
Let us denote by $\calQ_1$ the event that $\gamma_\cdot$ eventually settles in the first quadrant; notice that $P(\calQ_1)=1/4$. We then have, for any $\epsilon>0$,
\[
\lim_{t\to\infty} P\bigl( \sup_{0\leq s\leq t} \bigl\| \tfrac1t \widehat\gamma_s - \tfrac1t \gamma_s \bigr\|_{\scriptscriptstyle 2} \geq \epsilon \bgiven \calQ_1  \bigr) = 0.
\]
\end{theorem}
\begin{proof}
By Proposition~\ref{Prop1}, $\gamma_\cdot$ and $\widehat\gamma_\cdot$ couple almost surely in finite time. Since, on $\calQ_1$, the distance between the two processes (at the same time) remains constant after coupling, it follows that, almost surely, $\sup_{t\geq 0} \| \widehat\gamma_t - \gamma_t \|_{\scriptscriptstyle 1} < \infty$.
\end{proof}

\section{The scaling limit: Proof of Theorem~\ref{thm_main}}\label{sec_scalinglimit}
\begin{figure}
\begin{center} 
\scalebox{0.65}{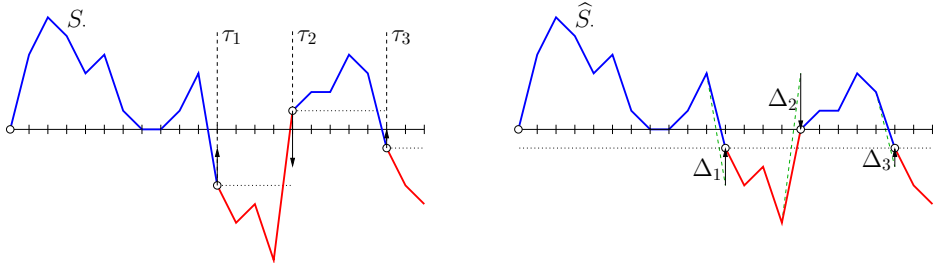}
\end{center}
\caption{The construction of the process $\hat S_\cdot$. }
\label{FigS_chapeau}
\end{figure}

Theorem~\ref{thm_main} is a consequence of Theorem~\ref{thm_cornervsprudent} and the following result.
\begin{theorem}\label{thm_slimit-corner}
On a suitably enlarged probability space, one can construct simultaneously a realization of $\widehat{\gamma}_\cdot$ and of a Brownian motion $W_\cdot$ such that, for any $\epsilon>0$,
\[
\lim_{t\to\infty} P \bigl(\sup_{0\leq s\leq t} \bigl\|\tfrac1t \widehat\gamma_s - Z^{1,1}_{s/t} \bigr\|_{\scriptscriptstyle 2} \geq \epsilon \bigr) = 0,
\]
where $Z^{1,1}_\cdot$ was defined in \eqref{processusZ}.
\end{theorem}
The rest of this section is devoted to the proof of Theorem~\ref{thm_slimit-corner}.
We start by establishing a suitable coupling between the corner process and the effective random walk $S_\cdot$ introduced in Subsection~\ref{ssec-effective}. To this end, we construct a process $\hat S_\cdot$ associated to the effective random walk $S_\cdot$. This construction is illustrated on Figure \ref{FigS_chapeau}.

\medskip
Consider the following alternating ladder times for $S_\cdot$: $\tau_0:=0$ and, for $k\in\bbZ_{\geq 0}$,
\begin{gather*}
\tau_{2k+1} := \inf\setof{n>\tau_{2k}}{S_n<S_{\tau_{2k}}},\\
\tau_{2k+2} := \inf\setof{n>\tau_{2k+1}}{S_n>S_{\tau_{2k+1}}}.
\end{gather*}
Observe that $\tau_{k+1}-\tau_k$ has the same distribution as $\eta_\infty$.
Define the overshoots as $\Delta_0:=0$, and for all $k\in\bbZ_{\geq 0}$,
\begin{align*}
\Delta_{2k+1} := -1 -(S_{\tau_{2k+1}}-S_{\tau_{2k}}),\\
\Delta_{2k+2} := +1 -(S_{\tau_{2k+2}}-S_{\tau_{2k+1}}).
\end{align*}

Thanks to the memoryless property of the geometric distribution, we have that, for all $k$, 
$(-1)^{k+1}\Delta_k$ is a non-negative random variable, with geometric distribution of parameter $1/2$. 
We now define
\begin{equation}\label{eq:hatSn}
\widehat S_n := S_n + \sum_{j\geq 0} \Delta_j \mathbf{1}_{\{\tau_j\leq n\}}. 
\end{equation}

The trajectory of the corner model can be constructed from the process $\widehat S_\cdot$, by decomposing the trajectory of the latter into two types of ``excursions'': those starting at $0$ and ending at $-1$, and those starting at $-1$ and ending at $0$. Applying the inverse of the procedure described in Fig.~\ref{Fig2} to an excursion of the first type, we construct a vertical excursion of $\widehat\gamma_\cdot$; similarly, applying the same procedure to an excursion of the second type, we construct a horizontal excursion of $\widehat\gamma_\cdot$. The fact that these have the proper distribution follows from Lemma~\ref{lem_effective}.

The next result shows that the sup-norm between the two processes $S_\cdot$ and $\widehat S_\cdot$ up to any fixed time never gets too large.
\begin{lemma}\label{lem:SwidehatS}
For any $\delta>0$,
\[
\lim_{n\to\infty} P\bigl(\max_{0\leq k\leq n} |S_k - \widehat S_k| \geq n^{1/4+\delta}\bigr) = 0.
\]
\end{lemma}
\begin{proof}
By \eqref{eq:hatSn}, we have $S_k - \widehat S_k=\tilde S_{N_{\rm o}(k)}-\Delta_{k+1}/2$, where 
$\tilde S_N:=\sum_{j=0}^{N} \tilde\Delta_i$, $\tilde \Delta_i:=(\Delta_{i}+\Delta_{i+1})/2$ and
$N_{\rm o}(k) := \max \setof{j\geq 0}{\tau_j \leq k}$. Observe that the random variables $\tilde\Delta_i$ are i.i.d.\ and symmetric; in particular, $\tilde S_N$ is a martingale. 
If $L:=n^{1/4+\delta}$,
\begin{align*}
\lim_{n\to\infty} P\bigl(\max_{0\leq k\leq n} |S_k - \widehat S_k| \geq n^{1/4+\delta}\bigr) &\leq
\lim_{n\to\infty} P\bigl(\max_{k\leq n}|\tilde S_{N_{\rm o}(k)}|\geq L/2\bigr)\\
&\leq  
\lim_{n\to\infty} P\bigl(\max_{j\leq N_{\rm o}(n)}|\tilde S_{j}|\geq L/2\bigr)\\
&\leq \lim_{n\to\infty} P\bigl(\max_{j\leq L^2}|\tilde S_{j}|\geq L/2\bigr)+\lim_{n\to\infty} P(N_{\rm o}(n) \geq L^2)\,.
\end{align*}
The first limit is zero by the Doob-Kolmogorov inequality, and so is the second by a standard renewal argument, since $P(\tau_1=n)\sim n^{-3/2}$.
\end{proof}
\begin{lemma}\label{lem:SS}
Fix $1/8>\delta>0$. For all $\epsilon>0$, we have
\begin{gather*}
\lim_{n\to\infty} P\Bigl( \max_{1\leq k \leq n} \bigl| \tfrac1n \sum_{i=1}^k \mathbf{1}_{\{\widehat S_i \geq 0\}} - \tfrac1n \sum_{i=1}^k \mathbf{1}_{\{S_i \geq n^{1/4+\delta}\}} \bigr| > \epsilon \Bigr) = 0,\\
\lim_{n\to\infty} P\Bigl( \max_{1\leq k \leq n} \bigl| \tfrac1n \sum_{i=1}^k \mathbf{1}_{\{\widehat S_i < 0\}} - \tfrac1n \sum_{i=1}^k \mathbf{1}_{\{S_i < -n^{1/4+\delta}\}} \bigr| > \epsilon \Bigr) = 0.
\end{gather*}
\end{lemma}
\begin{proof}
Fix $1/8>\delta>0$.
By the local central limit theorem, $E(\sum_{i=1}^n \mathbf{1}_{\{|S_i| < n^{1/4+\delta} \}}) \leq \Cl[c]{c40} n^{3/4+\delta}$. Therefore
\begin{equation}\label{eq:goodpoints}
\lim_{n\to\infty} P\Bigl( \sum_{i=1}^n \mathbf{1}_{\{|S_i| < n^{1/4+\delta} \}} \geq n^{3/4+2\delta} \Bigr) = 0.
\end{equation}
Observe now that the signs of $S_k$ and $\widehat S_k$ coincide when $|S_k| \geq n^{1/4+\delta}$, by Lemma~\ref{lem:SwidehatS}.
The conclusion follows, since
\begin{multline*}
P\Bigl( \max_{1\leq k \leq n} \bigl| \tfrac1n \sum_{i=1}^k \mathbf{1}_{\{\widehat S_i \geq 0\}} - \tfrac1n \sum_{i=1}^k \mathbf{1}_{\{S_i \geq n^{1/4+\delta}\}} \bigr| > \epsilon \Bigr)\\
\leq
P\Bigl( \max_{\epsilon n/2\leq k \leq n} \bigl| \tfrac1n \sum_{i=1}^k \mathbf{1}_{\{\widehat S_i \geq 0\}} - \tfrac1n \sum_{i=1}^k \mathbf{1}_{\{S_i \geq n^{1/4+\delta}\}} \bigr| > \epsilon/2 \Bigr).
\end{multline*}
\end{proof}
The next ingredient is the existence of a strong coupling between the effective random walk $S_\cdot$ and the Brownian motion $B_\cdot$~\cite{KMT76}: on a suitably enlarged probability space, we can construct both processes in such a way that, for $n$ large enough,
\[
P\bigl( \max_{k\leq n} |S_k - \sigma B_k| > n^{1/4} \bigr) \leq e^{-n^{1/4}},
\]
where $\sigma^2=2$ is the variance of $\xi_i$.
It thus follows from Lemma~\ref{lem:SwidehatS} that, for any $\delta>0$,
\begin{equation}\label{eq:SW}
\lim_{n\to\infty} P\bigl( \max_{k\leq n} |\widehat S_k - \sigma B_k| > n^{1/4 + \delta} \bigr) = 0.
\end{equation}
\begin{lemma}\label{lem:WS}
Fix $1/8>\delta>0$. For all $\epsilon>0$, we have
\begin{gather*}
\lim_{n\to\infty} P\Bigl( \max_{k\leq n} \bigl| \tfrac1n \int_0^{k} \mathbf{1}_{\{B_s \geq 0\}} \mathrm{d}s - \tfrac1n \sum_{i=1}^{k} \mathbf{1}_{\{S_i \geq n^{1/4+\delta}\}} \bigr| > \epsilon \Bigr) = 0,\\
\lim_{n\to\infty} P\Bigl( \max_{k\leq n} \bigl| \tfrac1n \int_0^{k} \mathbf{1}_{\{B_s < 0\}} \mathrm{d}s - \tfrac1n \sum_{i=1}^{k} \mathbf{1}_{\{S_i < -n^{1/4+\delta}\}} \bigr| > \epsilon \Bigr) = 0.
\end{gather*}
\end{lemma}
\begin{proof}
It follows from~\eqref{eq:goodpoints} that the number of times $k\leq n$ for which $|S_{k-1}|\wedge |S_{k}|\wedge |S_{k+1}| > n^{1/4+\delta}$ is at least $n-3n^{3/4+2\delta}$, with probability going to $1$ as $n\to\infty$. Of course, the same remains true if we impose additionally that $S_{k-1}$, $S_{k}$ and $S_{k+1}$ have the same sign. For such times $k$,
\[
P\Bigl( \mathbf{1}_{\{S_k \geq n^{1/4+\delta}\}} \neq \int_{k-1/2}^{k+1/2} \mathbf{1}_{\{B_s \geq 0\}} \mathrm{d}s \Bigr) \leq e^{-\Cl[c]{c41} n^{1/2 + 2\delta}},
\]
and the conclusion follows.
\end{proof}
Let us now introduce
\begin{gather*}
\theta^+(\widehat S_{[0,n]}) := \sum_{k=0}^n \mathbf{1}_{\{\widehat S_k \geq 0\}},\qquad \theta^-(\widehat S_{[0,n]}) := n - \theta^+(\widehat S_{[0,n]}),\\
\theta^+(B_{[0,t]}) := \int_0^t \mathbf{1}_{\{B_s \geq 0\}} \mathrm{d}s,\qquad \theta^-(B_{[0,t]}) := t - \theta^+(B_{[0,t]}).
\end{gather*}
Using these notations, we can deduce from Lemmas~\ref{lem:SS} and~\ref{lem:WS} that
\begin{equation}\label{eq:theta}
\lim_{n\to\infty} P\Bigl( \max_{k\leq n} \bigl| \tfrac1n \theta^+(B_{[0,k]}) - \tfrac1n \theta^+(\widehat S_{[0,k]}) \bigr| > \epsilon \Bigr) = 0,
\end{equation}
for all $\epsilon>0$, and similarly for $\theta^-(B_{[0,k]})$ and $\theta^-(\widehat S_{[0,k]})$.
Let us introduce, for $m\geq 0$,
\begin{gather*}
\Gamma_m := \theta^+(\widehat S_{[0,m]})\, \vec{e}_1 + \theta^-(\widehat S_{[0,m]})\, \vec{e}_2,\\
\sfZ_m := \theta^+(B_{[0,m]})\, \vec{e}_1 + \theta^-(B_{[0,m]})\, \vec{e}_2.
\end{gather*}
It follows from~\eqref{eq:theta} that, for all $\epsilon>0$,
\[
\lim_{n\to\infty} P\bigl( \sup_{0\leq m\leq n}  \bigl\| \tfrac1n\Gamma_m - \tfrac1n \sfZ_m \bigr\|_{\scriptscriptstyle 2} > \epsilon \bigr) = 0.
\]
Given $n\in\bbN$, we denote by $t(n):=\sum_{i=1}^n \bigl(1+|\widehat S_{i} - \widehat S_{i-1}|\bigr)$ the microscopic time such that the point $\widehat S_n$ is mapped on the point $\widehat\gamma_{t(n)}$ by the transformation described before Lemma~\ref{lem:SwidehatS}.
\begin{lemma}\label{lem:tube}
For any $\epsilon>0$,
\begin{equation}
\lim_{n\to\infty} P\bigl( \sup_{0\leq m\leq n} \bigl\| \tfrac1n \widehat\gamma_{t(m)} - \tfrac1n \Gamma_m \bigr\|_{\scriptscriptstyle 2} > \epsilon \bigr) = 0.
\end{equation}
\end{lemma}
\begin{proof}
This follows from Lemma~\ref{lem:SwidehatS} and the fact that, for any $\delta>0$, $\max_{0\leq k\leq n} |S_k| \leq n^{1/2+\delta}$, with probability going to $1$ as $n\to\infty$.
\end{proof}
It remains to relate more explicitly the real microscopic time $t(n)$ and the time $n$ of the effective random walk.
\begin{lemma}\label{lem:time}
For any $\epsilon >0$,
\[
\lim_{n\to\infty} P\bigl( \sup_{0\leq m\leq n} |t(m)-\tfrac73 m| > \epsilon n\bigr) = 0.
\]
\end{lemma}
\begin{proof}
As we have seen $\lim_{n\to\infty} P(N_{\rm o}(n)\geq n^{1/2+\delta}) = 0$, for all $\delta>0$. We thus deduce that
\begin{multline*}
\lim_{n\to\infty} P\bigl( \sup_{1\leq m \leq n}\bigl| \sum_{i=1}^m |\widehat S_{i} - \widehat S_{i-1}| - \sum_{i=1}^m |S_{i} - S_{i-1}| \bigr| > \epsilon n \bigr)\\
=
\lim_{n\to\infty} P\bigl( \sup_{1\leq m \leq n} \bigl| \sum_{i=1}^m |\widehat S_{i} - \widehat S_{i-1}| - \sum_{i=1}^m |S_{i} - S_{i-1}| \bigr| > \epsilon n, N_{\rm o}(n) < n^{1/2+\delta} \bigr).
\end{multline*}
Now, observing that
\[
\sup_{1\leq m\leq n}\bigl| \sum_{i=1}^m |\widehat S_{i} - \widehat S_{i-1}| - \sum_{i=1}^m |S_{i} - S_{i-1}| \bigr| \leq \sum_{i=1}^{N_{\rm o}(n)} |\Delta_i|,
\]
we deduce that
\begin{multline*}
\lim_{n\to\infty} P\bigl( \sup_{1\leq m \leq n}\bigl| \sum_{i=1}^m |\widehat S_{i} - \widehat S_{i-1}| - \sum_{i=1}^m |S_{i} - S_{i-1}| \bigr| > \epsilon n \bigr)\\
\leq
\lim_{n\to\infty} P\bigl(\sum_{i=1}^{n^{1/4+\delta}} |\Delta_i| > \epsilon n \bigr) = 0.
\end{multline*}
It then follows from the Doob-Kolmogorov inequality that, since $E(|S_i-S_{i-1}|)=E(|\xi_i|) = 4/3$,
\begin{multline*}
\lim_{n\to\infty} P\bigl( \sup_{1\leq m\leq n} |t(m)-\tfrac73 m| > \epsilon n\bigr)\\
\leq
\lim_{n\to\infty} P\bigl( \sup_{1\leq m\leq n} \bigl| \sum_{i=1}^m \bigl( 1+|S_{i} - S_{i-1}| - \tfrac73\bigr) \bigr| > \epsilon n\bigr) = 0.
\end{multline*}
\end{proof}
Actually, we rather need to express $n$ in terms of the real microscopic time $t$: $n(t):=\inf\setof{n\geq 0}{t(n)\geq t}$. However, since $t(n+1)-t(n) = |\xi_{n+1}|$, it follows from the previous lemma that $\lim_{t\to\infty} P(\sup_{0\leq s\leq t}|n(s)-\tfrac37 s| > \epsilon t) = 0$.

\medskip
To sum up, we have established that, for all $\epsilon>0$,
$$
\lim_{t\to\infty} P\bigl( \sup_{0\leq s\leq t} \bigl\| \tfrac1t\widehat\gamma_s - \tfrac1t \sfZ_{3s/7} \bigr\|_{\scriptscriptstyle 2} > \epsilon \bigr) = 0.
$$
The claim of Theorem~\ref{thm_slimit-corner} follows since, for $u\in[0,1]$,
\[
\tfrac1t \sfZ_{3ut/7} = \int_0^{3u/7} \bigl( \mathbf{1}_{\{W_s\geq 0\}} \vec{e}_1 + \mathbf{1}_{\{W_s< 0\}} \vec{e}_2 \bigr) \mathrm{d} s = Z_u^{1,1},
\]
where we have set $W_u := B_{u t}/\sqrt{t}$.
\section{Concluding remarks}
In the present work, we have focused on some of the most striking features of the scaling limit of the (kinetic) prudent walk. There remain however a number of open problems. We list a few of them here.
\begin{itemize}
\item It would be interesting to determine the scaling limit of the prudent walk on other lattices, e.g., triangular. Observe that the scaling limit we obtain reflects strongly the symmetries of the $\bbZ^2$ lattice, and is thus very likely to be different for other lattices. Nevertheless, numerical simulations indicate that similar scaling limits hold.
\item In this work, we relied heavily on properties specific to the 2-dimensional square lattice. There are two natural generalizations of the prudent walk 
in higher dimensions: (i) the natural extension (forbidding steps in directions where visited sites are present), (ii) forbidding steps such that the corresponding half-line intersect the bounding parallelepiped. Notice that both coincide in dimension $2$.

The second variant is easier, and it is likely to be in the scope of a suitable extension of our techniques. The first variant, however, is much more subtle (and interesting). Numerical simulations suggest the existence of an anomalous scaling exponent: If $\gamma_t$ is the location of the walk after $t$ steps, its norm appears to be of order $t^\alpha$ with $\alpha \simeq .75$ (but we see no reason to expect $\alpha$ to be equal to $3/4$).

%
\begin{figure}
\includegraphics[width=\textwidth]{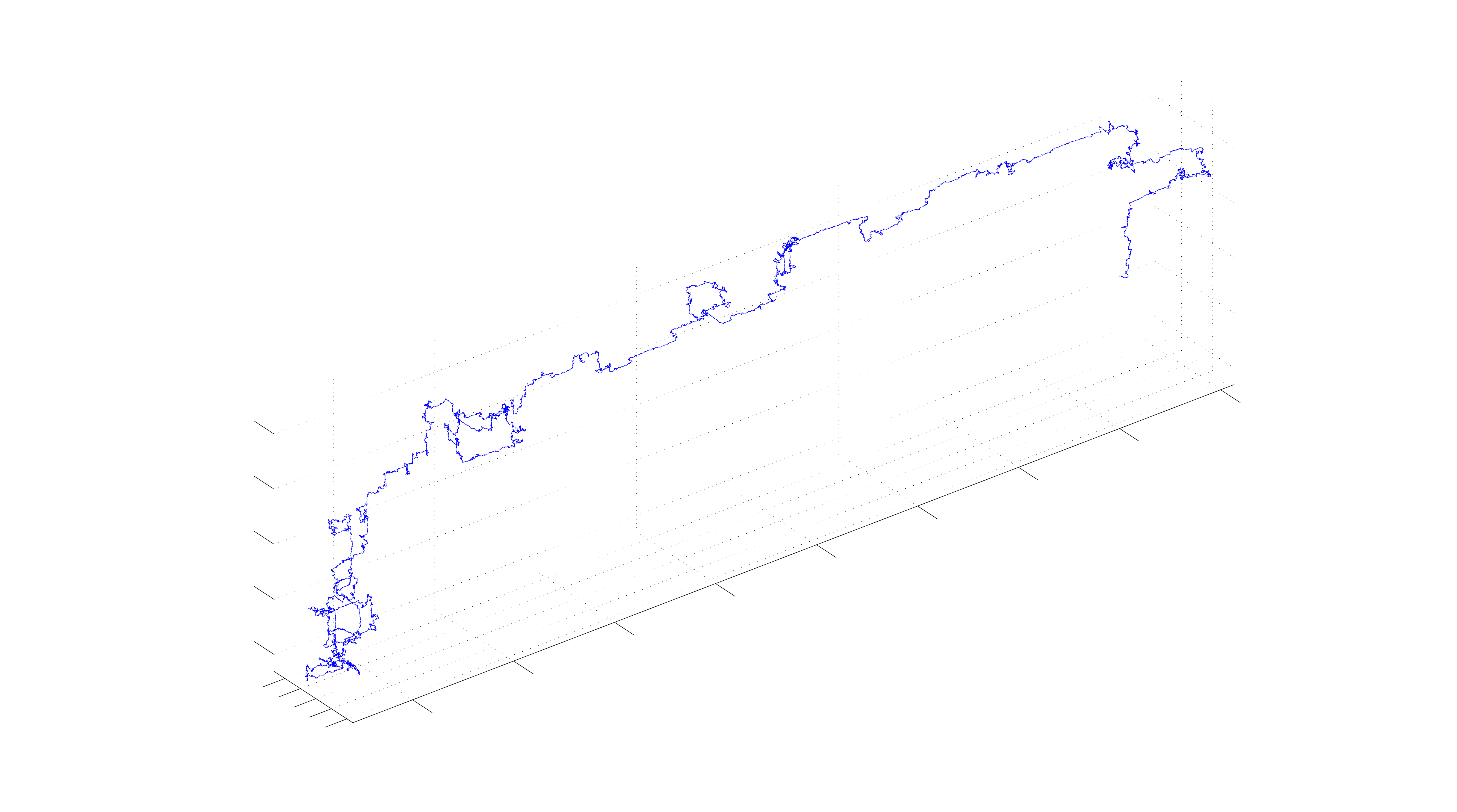}
\caption{A trajectory of a $3$-dimensional generalization of the prudent walk (first variant) with $6$ million steps.}
\label{Figure-3d}
\end{figure}
\item Derive the scaling limit of the uniform prudent walk. In particular, it was observed numerically~\cite{DGGJ2008} that the number of prudent walks has the same growth rate as that of the (explicitly computed) number of so-called 2-sided prudent walks (which are somewhat analogous to our corner model). This should not come as a surprise, since one expects the uniform prudent walk to visit all corners but one only finitely many times. It would be interesting to see whether it is possible to make sense of the fact that the uniform prudent walk is more strongly concentrated along (one of) the diagonal. This would permit to derive a proof of the previous result from the one of the kinetic model.
\end{itemize}

\bibliographystyle{plain}
\bibliography{BFV08}
\end{document}